\documentclass[11pt]{amsart}

\usepackage[utf8]{inputenc}
\usepackage{amssymb, amsmath}
\usepackage{xcolor}
\usepackage{graphicx}
\usepackage{hyperref}
\usepackage{geometry}
\usepackage{caption}
\usepackage{subcaption}

\usepackage[T1]{fontenc}
\definecolor{cmk2020}{HTML}{BE1E2D}
\hypersetup{
    colorlinks=true,
    citecolor=cmk2020,
    urlcolor=cmk2020
}
\newcommand{\pv}{\par\vspace{1ex}}

\theoremstyle{definition}

\theoremstyle{remark}

\numberwithin{equation}{section}

\newcommand{\R}{{\mathbb{R}}}

\newcommand{\N}{{\mathbb{N}}}

\newcommand{\F}{{\mathbb{F}}}

\begin{document}

\title{New model of non-Euclidean plane}

\author[Piotr B{\l}aszczyk, Anna Petiurenko]{Piotr B{\l}aszczyk, Anna Petiurenko}

\begin{abstract}We present a new model of a non-Euclidean plane, in which angles in a~triangle sum up to $\pi$. It is a subspace of the Cartesian plane over the field of hyperreal numbers $\R^*$. The model enables one to represent the negation of equivalent versions of the parallel axiom, such as the existence of the circumcircle of a triangle,  and Wallis' or Lagendre's axioms, as well as the difference between 
non-Euclidean and hyperbolic planes.

The model has unique educational advantages as expounding its crucial ideas requires only the basics of Cartesian geometry and non-Archimedean fields. 
\vspace*{0.5cm}

\centering
\noindent {\sc Keywords:} Parallel postulate, semi-Euclidean plane, hyperbolic geometry, hyperreal numbers\\
\end{abstract}

\maketitle
\thispagestyle{empty}
\setcounter{page}{1}


\tableofcontents

\newpage
\section{Introduction}

\begin{figure} [h!]
\centering
\includegraphics[scale=0.7]{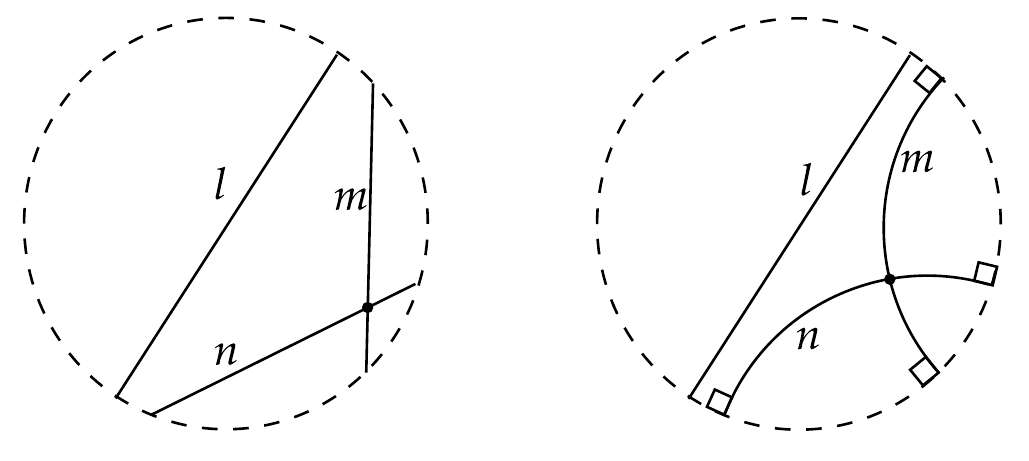}
\caption{Straight lines in Klein (left) and Poincar\'{e} (right) disk.} \label{poincare}
\end{figure}

Klein and Poincar\'{e} disks  are classical models of non-Euclidean geometry.  Both consist  of a~fixed  circle in the Euclidean plane, say $\Gamma$,  representing the  plane. In the Klein disk, chords of $\Gamma$ are straight lines; in the Poincar\'{e} disk, straight lines are diameters of $\Gamma$  or  arcs of circles orthogonal to $\Gamma$  (Fig. \ref{poincare}).

\begin{figure} [h!]
\centering
\includegraphics[scale=0.7]{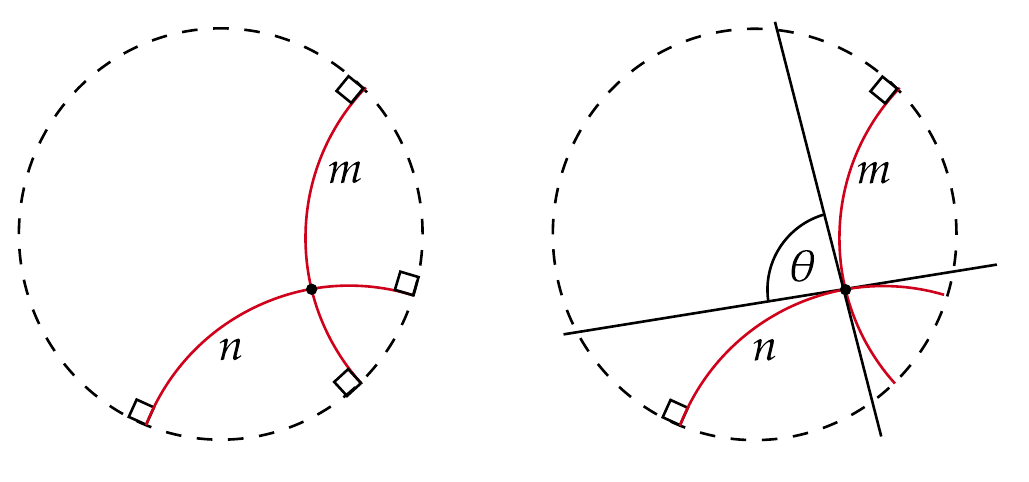}
\caption{Angles in the Poincar\'{e} disc.} \label{h_angle}
\end{figure}

\begin{figure} [!h]
\centering
\includegraphics[scale=0.7]{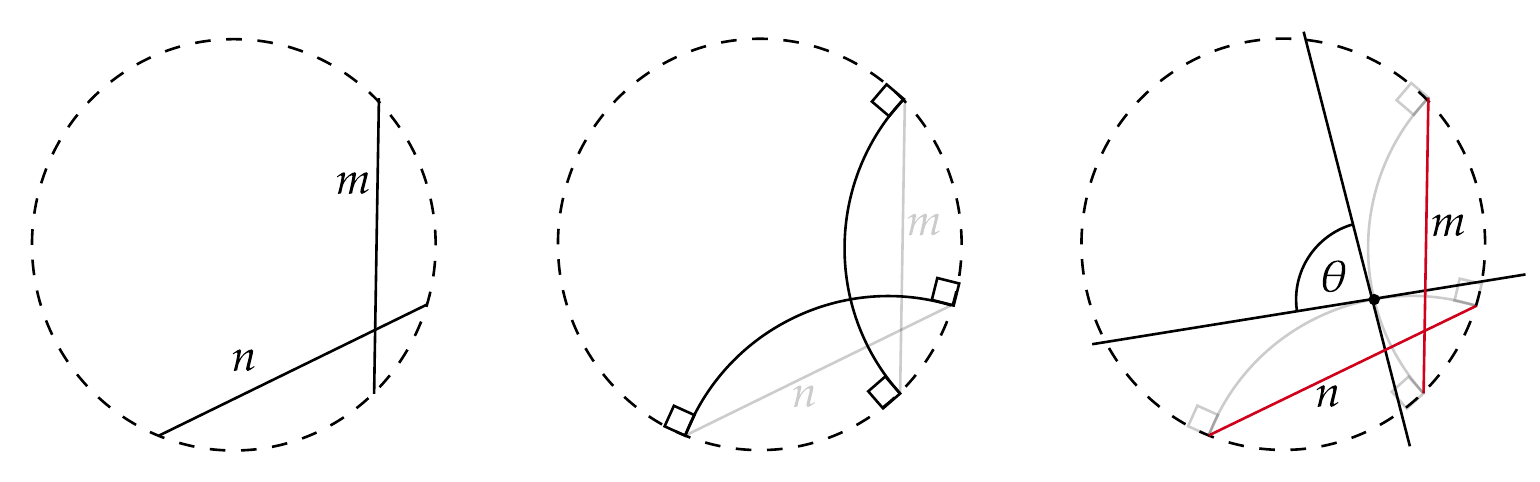}
\caption{Angles in the Klein disc.} \label{h_angle2}
\end{figure}

In the Poincar\'{e}  model, an angle between intersecting circles is the Euclidean angle between tangents to these lines drawn at their intersection point (Fig. \ref{h_angle}). In the Klein disc, an angle between intersecting straight lines is retrieved from the Poincar\'{e} model as presented in Fig. \ref{h_angle2}: for lines $n, m$, we draw circles orthogonal to $\Gamma$ and determine the angle between them.

Standard models of non-Euclidean planes, thus, involve a non-Euclidean representation of straight lines 
(Poincar\'{e}) or angles (Klein). We present a model in which both straight lines and angles are Euclidean as it takes the concept of a straight line and angle from the Cartesian plane over an ordered field.

A Euclidean plane is a Hilbert plane satisfying the parallel axiom and circle-circle intersection axiom. The Cartesian plane over an ordered field closed under the square root operation (i.e., over a Euclidean field)  provides a model of the Euclidean plane \cite[153]{ref_RH}.  Section 3 of our paper introduces a non-Archimedean field of hyperreals $(\R^*,+, \cdot,0,1,<)$. Since it is a Euclidean field, the Cartesian plane $\R^*\times \R^*$ is a model of Euclidean geometry.
We show that its subspace $\mathbb L\times \mathbb L$, where $(\mathbb L, +,\cdot, 0, 1)$ is the ring of limited hyperreal numbers,  is a semi-Euclidean plane, that is, it satisfies  Hilbert's axioms of the so-called absolute geometry, does not satisfy the parallel and Archimedean axioms, but angles in any triangle sum up to $\pi$ (i.e., two right angles). 

\section{Hilbert axioms for Euclidean and hyperbolic geometry}

Hilbert's \textit{Grundlagen der Geometrie}, from \cite{ref_DH99} to \cite{ref_DH72}  got eleven editions.
\cite{ref_RH} includes its modern version adjusted to educational practice.   Hilbert axioms, as presented therein, differ from the original only in applying  modern symbols.\footnote{\cite[597--602]{ref_MG} provides a~concise account of Hilbert axioms; see also \cite{ref_PB21}.}
 
 Hilbert grouped his axioms due to primitive concepts, point, straight line, plane, and the relation of betweenness, congruence of line segments, and angles. The paper will discuss three of them.

{\bf Archimedes' axiom}
 
 Given line segments $AB$ and $CD$, there is a~natural number $n$ such that $n$ 
copies of $AB$ added together will be greater than $CD$. 

{\bf Parallel axiom}

 For each point $A$ and each line $l$, there is at most one line containing $A$ that is parallel to $l$.
 
 \textbf{Hyperbolic axiom of parallels}

If $l$ is any line and $A$ a point not on the line, then there exist through $A$ two rays $AX^{\rightarrow}$ and $AY^{\rightarrow}$ which do not form one and the same line and  do not intersect the line $l$, while every ray emanating from $A$ that lies in the angle formed by $AX^{\rightarrow}$ and $AY^{\rightarrow}$  
does intersect $l$.\footnote{See  \cite[136]{ref_DH71}, \cite[259]{ref_MG}, \cite[374]{ref_RH}.}

\pv Absolute (neutral) geometry consists of axioms of incidence, betweenness, congruence of line segments and angles. One can also view it as a set of propositions  \textit{Elements} I.1--28.

Generally, parallel lines are not intersecting lines.  In the neutral geometry, the transportation of angles enables one to construct a line through $A$ that does not meet $l$. Hilbert's parallel axiom translates, thus, into the claim: There is exactly one line containing $A$ that is parallel to $l$.\footnote{\cite[75--77]{ref_PB21} discusses a relationship between Hilbert's and Euclid's parallel axioms.}  
In hyperbolic geometry, next to not intersecting lines, one considers limiting rays. 
Thus, Euclid's parallel axiom requires one line through $A$, not meeting $l$, its negation -- at least two lines through $A$ not meeting $l$, the hyperbolic axiom of parallels requires exactly two limiting rays originating in $A$, not meeting $l$. 
In what follows, we will present a non-Euclidean plane that is not hyperbolic.

\section{Semi-Euclidean plane}

In this section, we present a~model of a~semi-Euclidean plane, i.e., a~plane in which angles in a~triangle sum up to $\pi$ yet the parallel postulate fails. 
\cite[311]{ref_RH}, introduces that term, but the very idea originates in \cite[\S\,9]{ref_MD}. Dehn built such a~model owing to a~non-Archimedean Pythagorean field introduced in \cite[\S\,12]{ref_DH99};
yet, it was a~non-Euclidean field.\footnote{See also \cite[\S\,18]{ref_RH}. Example 18.4.3 expounds on Dehn's model.}  We employ the Euclidean field of hyperreal numbers.  To elaborate, let us start by introducing these numbers.

\subsection{The Cartesian plane over the field of hyperreal numbers }
An ordered field $(\F,+,\cdot,0,1,<)$ is a~ commutative field together with a~total order that is compatible with sums and products.
 In such a~field, one can define the following subsets of $\F$:
 \begin{enumerate}\itemsep 0mm
 \item[] $\mathbb L= \{x\in \F: (\exists n\in \N)(|x|<n)\}$,
 \item[] $\Psi= \{x\in \F: (\forall n\in \N)(|x|>n)\}$,
 \item[] $\Omega = \{x\in \F: (\forall n\in \N)(|x|<\tfrac 1n)\}$.
 \end{enumerate}

They are called limited, infinite, and infinitely small numbers, respectively.
 Here are some relationships helpful to pursue our arguments.
\begin{enumerate}\itemsep 0mm
\item[] $(\forall x, y\in\Omega)( x+y\in\Omega, x y\in\Omega)$,
\item[] $(\forall x, y\in\mathbb L)( x+y\in\mathbb L, x y\in\mathbb L)$,
\item[]  $(\forall x\in\Omega)( \forall y\in\mathbb L)( x y\in\Omega)$,
\item[] $(\forall x\neq 0)(x\in\Omega  \Leftrightarrow \ x^{-1}\in\Psi)$.
\end{enumerate}

We can also encode these rules by terms such as $\Omega+\Omega\subset \Omega$, or $\Omega\Omega\subset \Omega$, etc. 

To clarify our account, let us observe that the following equality $\Omega=\{0\}$ is a~version of the well-known    Archimedean axiom.

Since real numbers form the \textit{biggest} Archimedean field, every field extension of $(\R,+,\cdot,0,1,<)$
includes positive infinitesimals. The standard way of extending real numbers employs formal power series, below we sketch another method.

\begin{figure}[h!]
\centering
\includegraphics[scale=0.7]{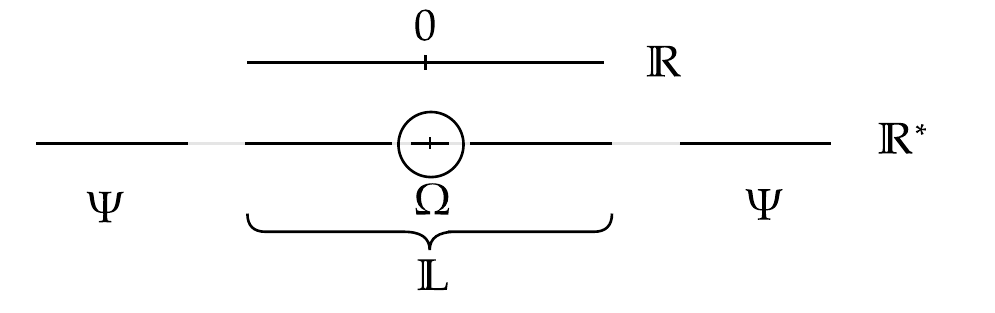}
\caption{The line of real numbers and its extension to hyperreals.} \label{figR*}
\end{figure}

 Let $\mathcal U$ be a~non-principal ultrafilter on $\N$.
The set of hyperreals is defined as a~reduced product $\R^*=\R^{\N}/{\mathcal U}$. Sums, products, and the order are introduced pointwise.
The field of hyperreals $(\R^*,+,\cdot,0,1,<)$
extends real numbers, hence, includes infinitesimals and infinite numbers \cite{ref_PB16}, \cite{ref_PB01}. 
Fig. \ref{figR*} represents in a~schematized way a~relationship between $\R$ and $\R^*$, as well as between $\mathbb L$, $\Psi$, and $\Omega$. 

To get an algebraic insight, let us note that limited numbers form an ordered ring,  while infinitesimals are its maximal ideal. Due to the so-called standard part theorem, one can show that the quotient ring is isomorphic to the field of real numbers. 
In consequence, the set $\mathbb L$ finds the following representation
\[\mathbb L=\bigcup_{r\in\R}r+\Omega.\] 

Although the below diagrams spot on the origin of the coordinate system, one can reiterate our arguments taking into account any point $(r,s)$, where $r,s\in\R$, and its infinitesimal neighborhood, i.e., $(r,s)+\Omega\times\Omega$, 
  in short $(r,s)+\Omega^2$.

In section \S\,3.3 below, we show how to introduce trigonometric functions and other counterparts of real maps, such as the square root. Thus, the field of hyperreals is a Euclidean field.

Due to the proposition  \cite[16.2]{ref_RH}, the Cartesian plane over the field of hyperreals is a~model of Euclidean plane, with straight lines and circles given by equations $ax+by+c=0$, $(x-a)^2+(y-b)^2=r^2$, 
where $a,b,c,r\in \R^*$ and due to the equation of a~straight line, parameters $a, b$ have to satisfy  condition $a^2+b^2\neq 0$; angles between straight lines are defined as in the  Cartesian plane over the field of real numbers. Specifically, on the plane $\R^*\times\R^*$, angles in triangles  sum up to $\pi$. Parallel lines are of the form $y=mx+b$ and $y=mx+c$, while a~perpendicular to the line $y=mx+b$ is given by the formula 
$y=-\frac{1}{m}x+d$. 

\begin{figure}[h!]
\centering
\includegraphics[scale=0.8]{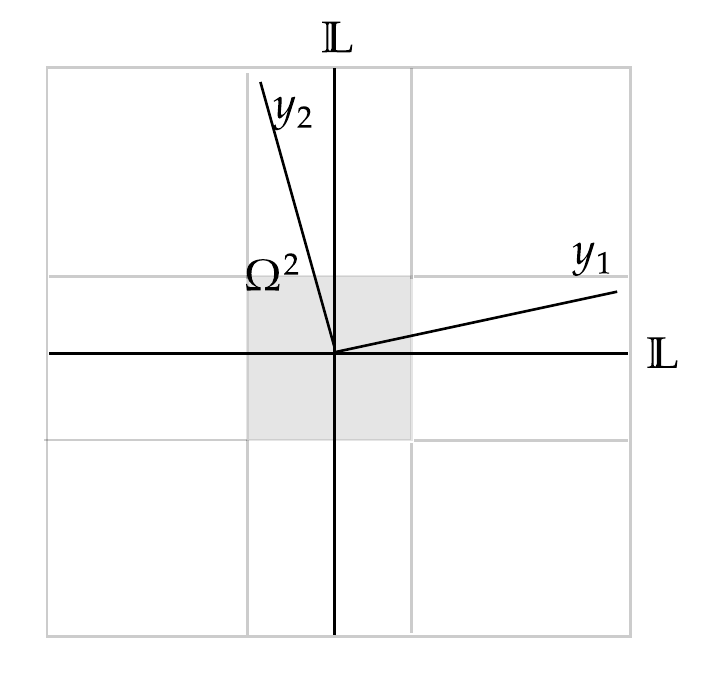}
\caption{Perpendicular lines with infinitesimal and infinitely large slopes.} \label{figLL}
\end{figure}

Now, let us take a~subspace $\mathbb L\times\mathbb L$ of the plane  $\R^*\times\R^*$. In that plane, circles are defined by analogous formula, namely 
$(x-a)^2+(y-b)^2=r^2$, where $a,b,c,r\in \mathbb L$, while  every line in $\mathbb L\times\mathbb L$ 
is of the form $l\cap \mathbb L\times\mathbb L$, where $l$ is a~line in $\R^*\times\R^*$. 
Since we want plane $\mathbb L\times\mathbb L$ include  lines such as $y_1=\varepsilon x$, where  
$\varepsilon\in\Omega$,  it also has to include the perpendicular $y_2=\frac{-1}{\varepsilon} x$, but 
$\frac{-1}{\varepsilon}\notin\mathbb L$. Formula 
$l\cap \mathbb L\times\mathbb L$, where $l=ax+by+c$ and $a, b, c \in \R^*$   guarantees the existence of the straight line $y_2$ in  $\mathbb L\times\mathbb L$. 
Finally, the interpretation of an angle is the same as in the model $\R^*\times \R^*$.

Explicit checking shows that the model characterized above satisfies Hilbert axioms of absolute plane geometry plus the circle-circle and line-circle axioms; see  \cite[\S 14--17]{ref_RH}.

\begin{figure}[h!]
\centering
\includegraphics[scale=0.7]{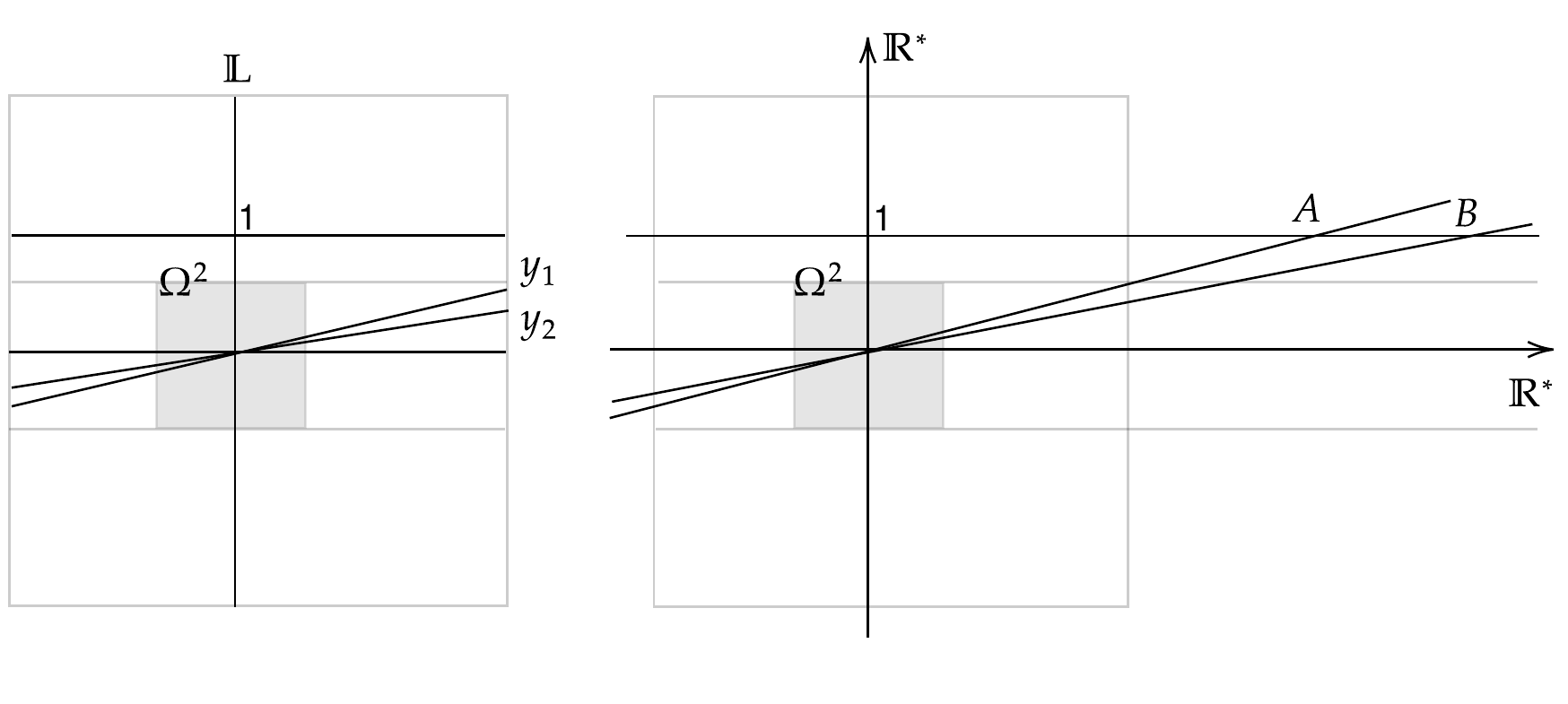}
\caption{Non-Euclidean plane $\mathbb L\times \mathbb L$ vs. Euclidean plane $\R^*\times \R^*$.} \label{figLLvsR*}
\end{figure}

With regard to parallel lines, let us consider the horizontal line $y\equiv 1$ and two specific lines through $(0,0)$,
namely $y_1=\varepsilon x, y_2=\delta x$, where $\varepsilon, \delta\in\Omega$; see Fig. \ref{figLLvsR*}.
Since
$\Omega\mathbb L\subset \Omega$, the following inclusions hold $y_1, y_2\subset \mathbb L\times \Omega$. In other words, values of maps $y_1(x), y_2(x)$ are infinitesimals, given that $x\in\mathbb L$.  The same obtains for any line of the form $y=\mu x$, with $\mu\in\Omega$. Since there are infinitely many infinitesimals, 
there are infinitely many lines through $(0,0)$ not intersecting the horizontal line $y\equiv 1$.

\begin{figure}[h!]
\centering
\includegraphics[scale=0.7]{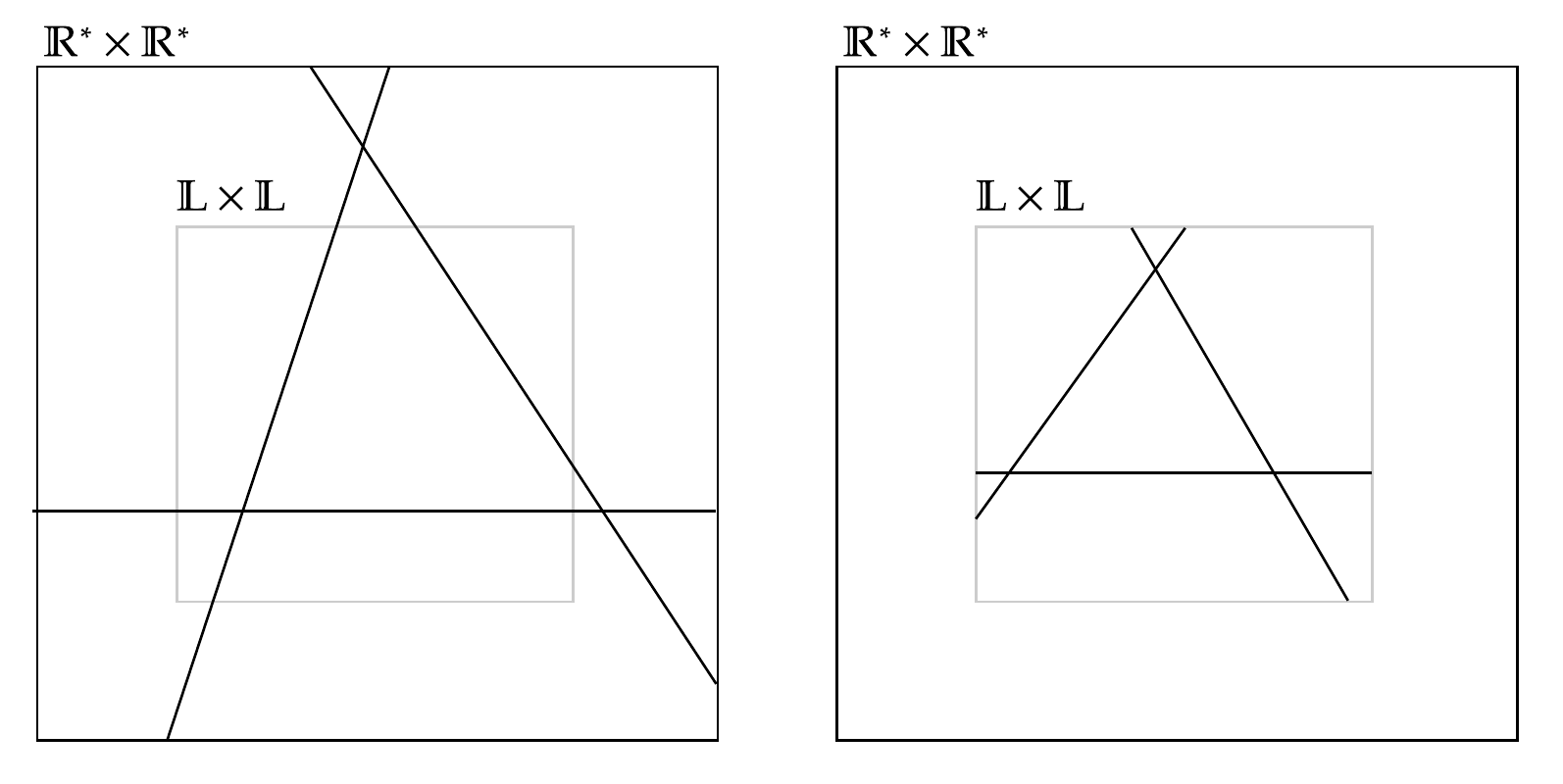}
\caption{Triangles in Euclidean plane $\R^*\times \R^*$ and its subspace $\mathbb L\times \mathbb L$.} \label{figT}
\end{figure}

Since every triangle in $\mathbb L\times \mathbb L$ is a~triangle in $\R^*\times \R^*$, it follows that angles in a~triangle  on the plane $\mathbb L\times \mathbb L$  sum up to $\pi$ (Fig. \ref{figT}).

\subsection{Circumcircle of a triangle, Wallis' and Lagendre's axioms}

In our model, one can easily find counterexamples for various versions of the parallel axiom.
Below, we discuss three of them.

      In the \textit{Elements}, Book IV, proposition  5, Euclid requires to circumscribe a circle about a given triangle. To this end, he   assumes that perpendicular bisectors of two sides of a triangle meet. Indeed, the intersection point is the center of the circle circumscribing a triangle.   The existence of such a circumcircle  is an equivalent version of the parallel axiom. 
 Below we show it does not hold in the plane  $\mathbb{L} \times \mathbb {L}$.

\begin{figure}[h!]
 \centerline{\includegraphics[scale=0.7]{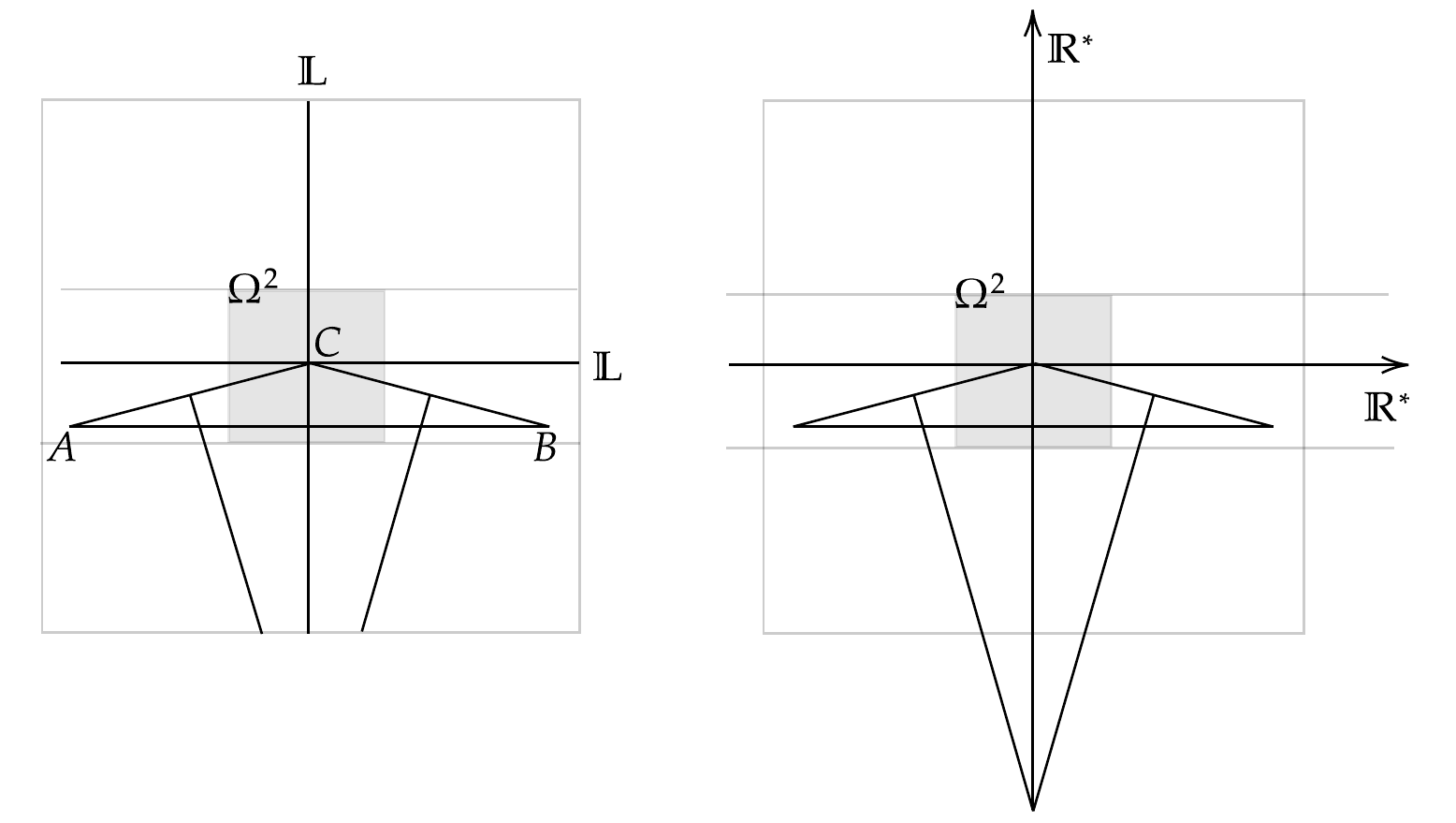}}
 \caption{Triangle in $\mathbb L\times\mathbb L$ with no circumcircle (left) and its counterpart in $\mathbb R^*\times\mathbb R^*$ (right).} \label{vsIV5}
 \end{figure}

 Let us take the line $y\equiv -\varepsilon$ and points $A=(-1, -\varepsilon)$, $B=(1, -\varepsilon)$ on it. A line   through the point $C=(0, 0)$ and $A$  has the equation $y=\varepsilon x$.
 Perpendicular bisectors of the sides $AB$ and $AC$ have the equations $x=0$  and 
 $2x +2\varepsilon y +\varepsilon^2 +1=0$, respectively. They meet at the point 
 $\left (0, -\varepsilon -\frac{1}{\varepsilon}\right)$. However,  $\left (0, -\varepsilon -\frac{1}{\varepsilon}\right)\notin \mathbb{L} \times \mathbb {L}$. Fig. \ref{vsIV5} depicts three perpendiculars to the sides of the triangle $ABC$.
 
 The so-called Wallis axiom reads: Given any triangle $\triangle ABC$ and given any segment $DE$, there exists a triangle  
 $\triangle DEF$ having $DE$ as one of his sides such that $\triangle ABC\sim \triangle DEF$ \cite[216]{ref_MG}.
 
 \begin{figure}[!h]
          \begin{subfigure}[b]{0.47\textwidth}
          \centering
        {\includegraphics[scale=0.7]{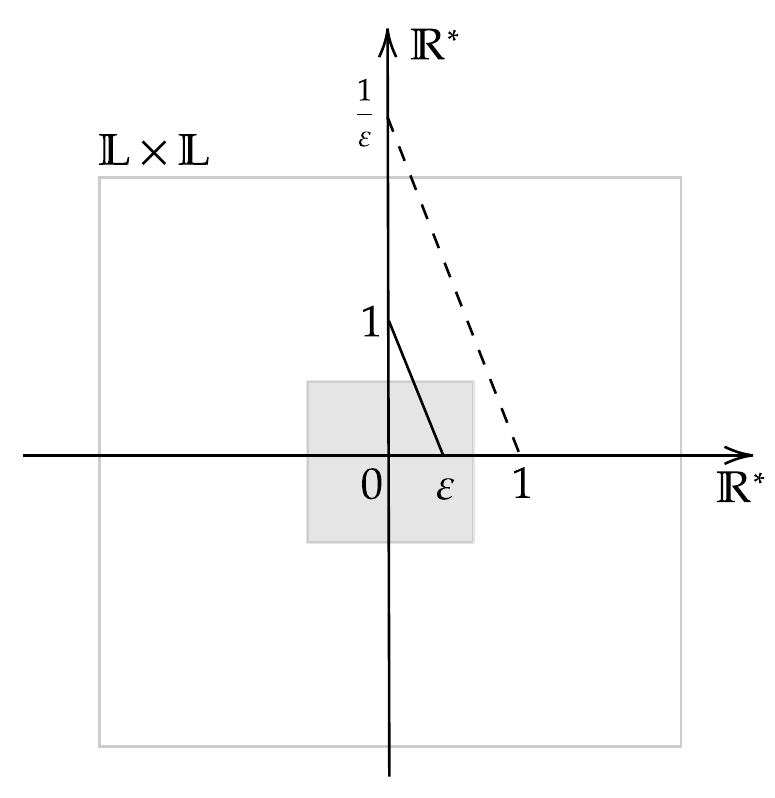}}
 \caption{Wallis axiom.} \label{WA}
      \end{subfigure}
     \hfill
     \begin{subfigure}[b]{0.47\textwidth}
     \centering
         {\includegraphics[scale=0.7]{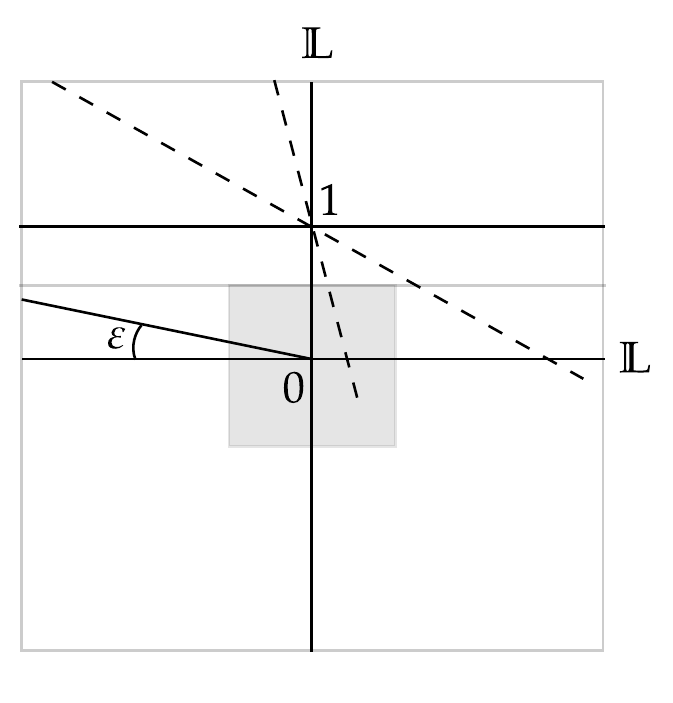}}
 \caption{Legender's Axiom.} \label{LA}
     \end{subfigure}
    \caption{}
\end{figure}
 
 It is yet another equivalent version of the parallel axiom \cite[216--217]{ref_MG}.

Let $ABC$ be a triangle with vertices $O$, $(\varepsilon,0)$, $(1,0)$, and take $DE=[0,1]$. Then, the triangle 
$DEF$ equiangular with  the triangle $ABC$ has the vertex $F=(0,\tfrac{1}{\varepsilon})$ outside the plane $\mathbb L\times\mathbb L$.

Finally, let us take Legender's axiom that reads: For any angle and any point $A$ in the interior of that angle, there exists a line through $A$ and not through that angle vertex which intersects both sides of the angle \cite[223]{ref_MG}.

 Lagendre's axiom is equivalent to the parallel axiom, given the Archimedean axiom \cite[324]{ref_RH}.
 
 To show it fails in our model, take $A=(0,1)$, and an angle with arms made by rays $y_1\equiv 0$ and  \mbox{$y_2=\tan(\pi-\varepsilon)x$} emanating from $O$. Now, any line cutting $y_1$ does not meet $y_2$. Moreover, the straight line $y_3\equiv 1$ does not meet any arm of the angle.
\subsection{Trigonometry}

Let $f$ be a~real map, i.e. $f:\R\mapsto \R$. Its extension to   a~map on $\R^*,$
  $f^*:\R^*\mapsto \R^*$, is defined by
\begin{equation}\label{nsa2}f^*([(r_n)]) =[(f(r_1),f(r_2),...)].\end{equation}

If $r\in \R$, then $f^*(r)=[(f(r),f(r),...)]$.
  Since we identify real number $r$ with hyperreal $[(r,r,...)]$,  the equality $f^*([(r,r,...)])=f(r)$ obtains, meaning, $f^*$ extends $f$, $f^*_{|\R}=f$.

Putting $f=\sqrt{\ \ }$ in definition (\ref{nsa2}), we get
$$\sqrt{[(r_n)]}^*= [(\sqrt{r_n})]=[(\sqrt{r_1}, \sqrt{r_2}, ...)],\ \ \mbox{for}\ \ [(r_j)]>0.$$


Similarity, under the definition (\ref{nsa2}), we have
\[\sin^*[(r_n)]=[(\sin r_1,\sin r_2,...)], \ \ \cos^*[(r_n)]=[(\cos r_1,\cos r_2,...)].\]

Since for every $n$ the identity $\sin^2 r_n+\cos^2 r_n=1$  holds, we have 
\[(\sin^*x)^2+(\cos^* x)^2=1.\]

Similarly, every trigonometric identity translates into an identity involving the maps $\sin^*$  and  $\cos^*$, $\tan^*$ and $\cot^*$. Interpreting angles and rotations in $\R^*\times \R^*$, one can apply these $^*$maps. 
Yet, that way of
 introducing trigonometry is hardly popular
among geometry scholars.

While dealing with angles in a Cartesian plane over a Pythagorean field, i.e., closed under the operation 
$x\mapsto \sqrt{1+x^2}$, Hilbert's solution was the following \cite[\S 9]{ref_DH99}. Let points 
$O, E, C$ have coordinates $(0,0)$, $(1,0)$,  $(c_x,c_y)$. Then  the rotation of the angle $\angle EOC$ about $O$ turns point  $(x,y)$ into $(x',y')$ defined as follows
\[x'=\frac{c_x}{\sqrt{c_x^2+c_y^2}}x-\frac{c_y}{\sqrt{c_x^2+c_y^2}}y,\]
\[y'=\frac{c_y}{\sqrt{c_x^2+c_y^2}}x+\frac{c_x}{\sqrt{c_x^2+c_y^2}}y.\]
  
Thus, terms $\frac{c_x}{\sqrt{c_x^2+c_y^2}}x$ and $\frac{c_y}{\sqrt{c_x^2+c_y^2}}y$ are to mimic the cosine and sine of the angle $\angle COE$, nonetheless, they make sense in any Pythagorean field. 

Considering angles and rotations in a Cartesian plane, Hartshorne mimics the real function $tan$.\footnote{See \cite[\S 16--17]{ref_RH}.}  Greenberg excludes the Poincare length in a non-Archimedean, Euclidean field.\footnote{See \cite[320]{ref_MG}.} 
Max Dehn introduced semi-Euclidean planes over non-Archimedean Pythagorean fields. And indeed, for a long time, non-Archimedean, Euclidean fields have not been known, and it could be a reason why Hilbert did not introduce the circle-circle axiom into his system of axioms. The only non-Archimedean fields known at that time were Pythagorean fields.  

 \section{Limiting parallel rays}

In this section we show that the plane $\mathbb L\times \mathbb L$ is not hyperbolic by considering limiting parallel rays. 

 \begin{figure}[h!]
 \centerline{\includegraphics[scale=0.7]{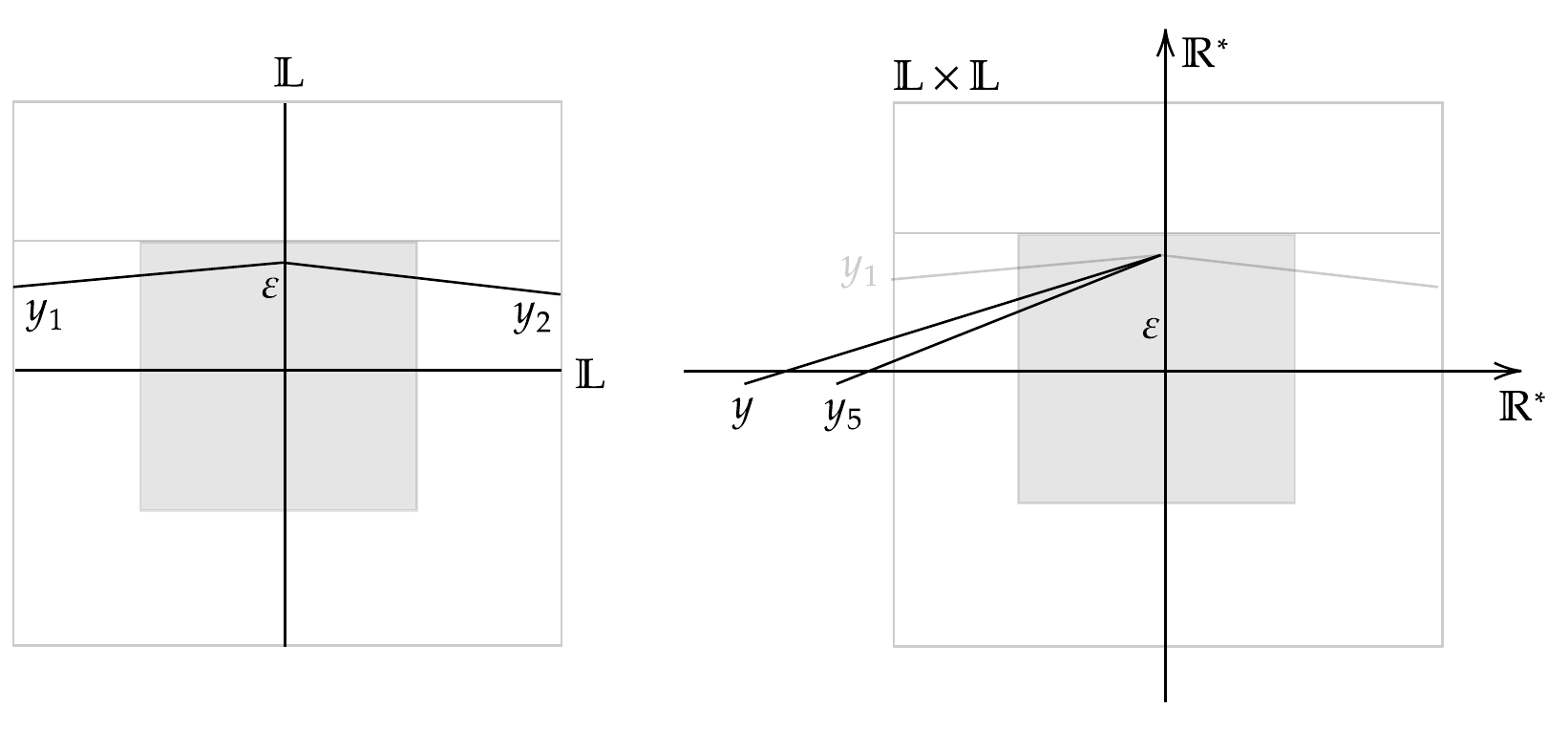}}
 \caption{Limiting parallel rays}
 \end{figure}

Take the point $A=(0,\varepsilon)$, for some positive infinitesimal $\varepsilon$.  Rays $y_1=\varepsilon^2x+\varepsilon$ and $y_2=-\varepsilon^2x+\varepsilon$ contain $A$ and do not meet $y_3\equiv 0$. 
Suppose 
$y=\delta x+\varepsilon$ is the limiting ray, that is, such that every line inside the angle with arms $y$ and $AO$ 
meets $y_3$. Slope $\delta$ has to satisfy the condition $\tfrac {-\varepsilon}{\delta}\in \Psi$. However, 
$y_5=\sqrt \delta x+\varepsilon$ lies inside the angle with arms $y$ and $AO$, contains $A$ and does not meet $y_3$. Therefore, $\mathbb L\times \mathbb L$ is not a hyperbolic plane.
 
Usually, non-Euclidean and not hyperbolic planes are modeled through some tricks concerning the algebraic properties of a field. 
It could be, for example, a Poincar\'{e} disc over a Pythagorean field that is not Euclidean.  Our plane $\mathbb L\times\mathbb L$ is not hyperbolic because it is 
semi-Euclidean. Nevertheless, one can also observe that both Euclidean and hyperbolic parallel axioms fail because straight lines in that plane are too `short'.


\end{document}